\documentclass[12pt]{amsart}
\usepackage{amssymb,latexsym}

\newtheorem{theorem}{Theorem}[section]

\theoremstyle{definition}

\def\ZZ{\mathbb{Z}}

\def\QQ{\mathbb{Q}}

\def\Acal{\mathcal{A}}

\renewcommand{\eqref}[1]{{\rm (\ref{#1})}}
\newcommand{\rem}[1]{\left\langle#1\right\rangle}

\begin{document}

\title[Semicanonical basis for cluster algebra of type $A_1^{(1)}$]
{Semicanonical basis generators of the cluster algebra of type
$A_1^{(1)}$}

\author{Andrei Zelevinsky}
\address{\noindent Department of Mathematics, Northeastern University,
 Boston, MA 02115}
\email{andrei@neu.edu}


\date{June 29, 2006}

 \thanks{Research supported by NSF (DMS) grant \# 0500534
and by a Humboldt Research Award.}

\maketitle

\section{Introduction}

The (coefficient-free) cluster algebra $\Acal$ of type $A_1^{(1)}$
is a subring of the field $\QQ(x_1, x_2)$ generated by the
elements $x_m$ for $m \in \ZZ$ satisfying the recurrence relations
\begin{equation}
\label{eq:x-recursion}
x_{m-1} x_{m+1} = x_m^2+1 \quad (m \in \ZZ)\ .
\end{equation}
This is the simplest cluster algebra of infinite type; it
was studied in detail in \cite{calzel,sherzel}.
Besides the generators $x_m$ (called \emph{cluster variables}),
$\Acal$ contains another important family of elements $s_0, s_1, \dots$
defined recursively by
\begin{equation}
\label{eq:s-recursion}
s_0 = 1, \,\, s_1 = x_0 x_3 - x_1 x_2, \,\
s_n = s_1 s_{n-1} - s_{n-2} \quad (n \geq 2).
\end{equation}
As shown in \cite{calzel,sherzel}, the elements $s_1, s_2, \dots$
together with the \emph{cluster monomials} $x_m^p x_{m+1}^q$ for
all $m \in \ZZ$ and $p, q \geq 0$, form a $\ZZ$-basis of $\Acal$
referred to as the \emph{semicanonical basis}.

As a special case of the \emph{Laurent phenomenon} established in
\cite{ca1}, $\Acal$ is contained in the Laurent polynomial
ring $\ZZ[x_1^{\pm 1}, x_2^{\pm 1}]$.
In particular, all $x_m$ and $s_n$ can be expressed as integer Laurent
polynomials in $x_1$ and $x_2$.
These Laurent polynomials were explicitly computed in \cite{calzel}
using their geometric interpretation due to P.~Caldero and F.~Chapoton
\cite{caldchap}.
As a by-product, there was given a combinatorial interpretation of these Laurent
polynomials, which can be easily seen to be equivalent to the one
previously obtained by G.~Musiker and J.~Propp
\cite{musikerpropp}.

The purpose of this note is to give short, self-contained and
completely elementary proofs of the combinatorial interpretation
and closed formulas for the Laurent polynomial expressions
of the elements $x_m$ and $s_n$.

\section{Results}

We start by giving an explicit combinatorial expression for each $x_m$ and
$s_n$, in particular proving that they are Laurent polynomials in
$x_1$ and $x_2$ with positive integer coefficients.
By an obvious symmetry of relations \eqref{eq:x-recursion}, each element $x_m$ is
obtained from $x_{3-m}$ by the automorphism of the ambient field
$\QQ(x_1, x_2)$ interchanging $x_1$ and $x_2$.
Thus, we restrict our attention to the elements $x_{n+3}$ for $n \geq 0$.

Following  \cite[Remark~5.7]{calzel} and \cite[Example~2.15]{yga},
we introduce a family of \emph{Fibonacci
polynomials} $F(w_1, \dots, w_N)$ given by
\begin{equation}
\label{eq:Fw}
F(w_1, \dots, w_N) = \sum_D \prod_{k \in D} w_k,
\end{equation}
where $D$ runs over all \emph{totally disconnected} subsets of $\{1, \dots, N\}$,
i.e., those containing no two consecutive integers.
In particular, we have
$$F(\emptyset) = 1, \,\, F(w_1) = w_1 + 1, \,\,
F(w_1, w_2) = w_1 + w_2 + 1.$$
We also set
\begin{equation}
\label{eq:fN}
f_N = x_1^{-\lfloor \frac{N+1}{2}\rfloor}
x_2^{-\lfloor \frac{N}{2}\rfloor} F(w_1, \dots, w_{N})|_{w_k = x_{\rem{k+1}}^2},
\end{equation}
where $\rem{k}$ stands for the element of $\{1,2\}$ congruent to $k$ modulo~$2$.
In view of \eqref{eq:Fw}, each $f_N$ is a Laurent polynomial in
$x_1$ and $x_2$ with positive integer coefficients.
In particular, an easy check shows that
\begin{equation}
\label{eq:f012}
f_0 = 1, \quad f_1 = \frac{x_2^2 + 1}{x_1} = x_3, \quad
f_2 = \frac{x_1^2+x_2^2+1}{x_1x_2}  = s_1.
\end{equation}

\begin{theorem}\cite[Formula~(5.16)]{calzel}
\label{th:xs-thru-f}
For every $n \geq 0$, we have
\begin{equation}
\label{eq:xs-thru-f}
s_n = f_{2n}, \quad
x_{n+3} = f_{2n+1}.
\end{equation}
In particular, all $x_m$ and $s_n$ are Laurent polynomials in
$x_1$ and $x_2$ with positive integer coefficients.
\end{theorem}

Using the proof of Theorem~\ref{th:xs-thru-f}, we derive the
explicit formulas for the elements~$x_m$ and $s_n$.

\begin{theorem}\cite[Theorems~4.1, 5.2]{calzel}
\label{th:xs-formula}
For every $n \geq 0$, we have
\begin{eqnarray}
\label{eq:xn-formula}
&x_{n+3} & =  x_1^{-n-1} x_2^{-n} (x_2^{2(n+1)} +
\sum_{q + r \leq n} {n-r \choose q}{n+1-q \choose r} x_1^{2q}
x_2^{2r});\\
\label{eq:sn-formula}
&s_{n} & =  x_1^{-n} x_2^{-n} \sum_{q+r \leq n}
{n-r \choose q}{n-q \choose r} x_1^{2q} x_2^{2r}.
\end{eqnarray}
\end{theorem}

\section{Proof of Theorem~\ref{th:xs-thru-f}}

In view of \eqref{eq:Fw}, the Fibonacci polynomials satisfy the
recursion
\begin{equation}
\label{eq:Fw-recursion}
F(w_1, \dots, w_N) =  F(w_1, \dots, w_{N-1}) + w_N F(w_1, \dots, w_{N-2})
\quad (N \geq 2).
\end{equation}
Substituting this into \eqref{eq:fN} and clearing the
denominators, we obtain
\begin{equation}
\label{eq:fN-recursion}
x_{\rem{N}} f_N =  f_{N-1} + x_{\rem{N-1}} f_{N-2}
\quad (N \geq 2).
\end{equation}
Thus, to prove \eqref{eq:xs-thru-f} by induction on~$n$, it
suffices to prove the following identities for all $n \geq 0$
(with the convention $s_{-1} = 0$):
\begin{eqnarray}
\label{eq:x1-xn}
&x_1 x_{n+3} &= s_{n} + x_2 x_{n+2};\\
\label{eq:x2-sn}
&x_2 s_n &= x_{n+2} + x_1 s_{n-1}.
\end{eqnarray}
We deduce \eqref{eq:x1-xn} and \eqref{eq:x2-sn} from \eqref{eq:s-recursion}
and its analogue established in \cite[formula~(5.13)]{sherzel}:
\begin{equation}
\label{eq:s1-xn}
x_{m+1} = s_1 x_{m} - x_{m-1} \quad (m \in \ZZ).
\end{equation}
(For the convenience of the reader, here is the proof of \eqref{eq:s1-xn}.
In view of \eqref{eq:s-recursion} and \eqref{eq:x-recursion}, we have
$$s_1 = \frac{x_1^2+x_2^2+1}{x_1x_2} = \frac{x_1+x_3}{x_2}
= \frac{x_0+x_2}{x_1}.$$
By the symmetry of the relations  \eqref{eq:x-recursion}, this
implies that $s_1 = (x_{m-1}+x_{m+1})/x_m$ for all $m \in \ZZ$,
proving \eqref{eq:s1-xn}.)

We prove \eqref{eq:x1-xn} and \eqref{eq:x2-sn} by induction on~$n$.
Since both equalities hold for $n = 0$  and $n=1$, we can assume
that they hold for all $n < p$ for some $p \geq 2$, and it
suffices to prove them for $n = p$.
Combining the inductive assumption with \eqref{eq:s-recursion}
and \eqref{eq:s1-xn}, we obtain
\begin{align*}
x_1 x_{p+3} &= x_1 (s_1 x_{p+2} - x_{p+1})\\
&= s_1(s_{p-1} + x_2 x_{p+1}) - (s_{p-2} + x_2 x_{p})\\
&= (s_1 s_{p-1} - s_{p-2}) + x_2(s_1 x_{p+1} - x_p)\\
&= s_p + x_2 x_{p+2},
\end{align*}
and
\begin{align*}
x_2 s_{p} &= x_2 (s_1 s_{p-1} - s_{p-2})\\
&= s_1(x_{p+1} + x_1 s_{p-2}) - (x_{p} + x_1 s_{p-3})\\
&= (s_1 x_{p+1} - x_{p}) + x_1(s_1 s_{p-2} - s_{p-3})\\
&= x_{p+2} + x_1 s_{p-1},
\end{align*}
finishing the proof of Theorem~\ref{th:xs-thru-f}.

\section{Proof of Theorem~\ref{th:xs-formula}}

Formulas \eqref{eq:xn-formula} and \eqref{eq:sn-formula}
follow from \eqref{eq:x1-xn} and \eqref{eq:x2-sn} by induction on~$n$.
Indeed, assuming that, for some $n \geq 1$, formulas \eqref{eq:xn-formula}
and \eqref{eq:sn-formula} hold for all the terms on the right hand
side of \eqref{eq:x1-xn} and \eqref{eq:x2-sn}, we obtain
\begin{eqnarray*}
x_{n+3} &=& x_1^{-1} (s_{n} + x_2 x_{n+2})\\
&=&x_1^{-n-1} x_2^{-n}(\sum_{q+r \leq n}
{n-r \choose q}{n-q \choose r} x_1^{2q} x_2^{2r}\\
&&+  (x_2^{2(n+1)} +
\sum_{q + r \leq n-1} {n-1-r \choose q}{n-q \choose r} x_1^{2q}
x_2^{2(r+1)}))\\
& =&  x_1^{-n-1} x_2^{-n} (x_2^{2(n+1)} +
\sum_{q + r \leq n} {n-r \choose q}({n-q \choose r} +
{n-q \choose r-1}) x_1^{2q} x_2^{2r})\\
& =&  x_1^{-n-1} x_2^{-n} (x_2^{2(n+1)} +
\sum_{q + r \leq n} {n-r \choose q}{n+1-q \choose r} x_1^{2q}
x_2^{2r}),
\end{eqnarray*}
and
\begin{eqnarray*}
s_{n} &=& x_2^{-1} (x_{n+2} + x_1 s_{n-1})\\
&=&x_1^{-n} x_2^{-n}(x_2^{2n} +
\sum_{q + r \leq n-1} {n-1-r \choose q}{n-q \choose r} x_1^{2q}
x_2^{2r}\\
&&+ \sum_{q+r \leq n-1}
{n-1-r \choose q}{n-1-q \choose r} x_1^{2(q+1)} x_2^{2r})\\
& =&  x_1^{-n} x_2^{-n} \sum_{q + r \leq n} ({n-1-r \choose q} +
{n-1-r \choose q-1}) {n-q \choose r} x_1^{2q} x_2^{2r}\\
& =& x_1^{-n} x_2^{-n} \sum_{q+r \leq n}
{n-r \choose q}{n-q \choose r} x_1^{2q} x_2^{2r},
\end{eqnarray*}
as desired.

\end{document}